\font\sixrm=cmr6
\font\eightrm=cmr8
\font\eighti=cmmi8

\font\tenbb=msbm10

\font\fourteenbf=cmbx12 at 14pt
\font\eightbf=cmbx8

\font\eighttt=cmtt8

\font\eightit=cmti8

\font\eightbfti=cmbxti10 at 8pt

\font\tensc=cmcsc10

\newfam\bbfam
\textfont\bbfam=\tenbb
\def\bb{\fam\bbfam\tenbb}

\newfam\scfam
\textfont\scfam=\tensc
\def\sc{\fam\scfam\tensc}

\def\eightpt{%
\def\rm{\eightrm}%
\def\it{\eightit}%
\def\tt{\eighttt}%
\def\bf{\eightbf}%
\def\bfti{\eightbfti}%
\textfont1=\eighti%
}


\def\NN{{\bb N}}

\def\CC{{\bb C}}
\def\Real{\mathop{\rm Re}}

\def\cF{{\cal F}}
\def\cG{{\cal G}}
\def\cV{{\cal V}}
\def\wh{\widehat}

\def\Un{{\bf 1}}
\def\hb{\hfil\break}


\parindent=0pt
\hsize= 6in
\hoffset= 0.25in

\tenrm
\baselineskip=22pt 
\parskip=0pt

\noindent{}

{\baselineskip=18pt\leftskip=1cm\rightskip=1cm

\vskip 1cm

\hfill{\tentt math.NT/0105120}\par

\vskip 1cm

{\fourteenbf
Sur certains espaces de Hilbert de fonctions enti\`eres, li\'es \`a la
transformation de Fourier et aux fonctions L de Dirichlet et de Riemann\par}

\vskip 0.5 in

\centerline{Jean--Fran\c{c}ois B{\eightrm URNOL}}
\centerline{Mai 2001}

\vskip 0.5 in
\baselineskip=16pt

Nous construisons dans un Espace de Sonine de fonctions enti\`eres un
sous-espace attach\'e \`a la fonction dz\^eta de Riemann et nous montrons que le
quotient contient des vecteurs intrins\`equement li\'es aux z\'eros non-triviaux
et \`a leurs multiplicit\'es \'eventuelles.

\vfill
{\leftskip=2cm
\baselineskip=12pt
Universit\'e de Nice~Sophia-Antipolis\hb
Laboratoire J.~A.~Dieudonn\'e, Math\'ematiques\hb
Parc Valrose\hb
F--06108 Nice cedex~02\hb
France\hb
\tt burnol@math.unice.fr\par}

\vfill

\eightrm
Cette Note a \'et\'e pr\'esent\'ee aux C. R. Acad. Sci., Paris, par le Professeur
Jean-Pierre K{\sixrm AHANE}. 

}

\eject

\noindent{}

\vskip1cm

\parindent = 5mm

\centerline{\fourteenbf English summary}

\bigskip

Let $K = L^2((0,\infty),dt)$. The Mellin transform of a function $f\in K$ is
defined, on the critical line $\Real(s) = 1/2$ and with the meaning of a
Plancherel isometry, through the formula $\wh{f}(s) = \int_0^\infty
f(t)\,t^{s-1}\,dt$.  Let $\cF_+$  be the cosine transform (formally
$\cF_+(f)(t)=2\int_0^\infty \cos(2\pi tu)f(u)du$.) Let $a$ and $b$ be two
strictly positive real numbers and let $K_{a,b}$ be the closed subspace of $K$
consisting of functions $f$ with support in $(a,\infty)$ and such that
$\cF_+(f)$ has its support in $(b,\infty)$. It is known classically that
$L^2((0,a])+\cF_+(L^2((0,b]))$ is a closed proper subspace of $K$ ({\it see\/}
[15], {\it sect. 2.9, p. 126-127\/}) hence that none of the $K_{a,b}$'s is
reduced to $\{0\}$. The functions
$\pi^{-(1-s)/2}\Gamma((1-s)/2)\wh{f}(s)$ on the critical line (for $f\in
K_{a,b}$) build (up to a trivial change of variable) a special case of a Sonine
Hilbert space as considered by De~Branges in [13], and in particular it consists of
{\it entire} functions. We prove elementarily:

\medskip

{\sc Th\'eor\`eme 1.\quad} {\it The Mellin transforms  of the
functions in $K_{a,b}$ are entire functions. They have at $s = 1+2j$, $j\in \NN$,
trivial zeros and these are their only common zeros (and are simple as such.) }

\medskip

We write $(\phi,\psi]$ for the {\it euclidean} scalar product $\int
\phi(u)\psi(u)du$. Let $I$ be the operator $I(f)(t) = 1/t\;f(1/t)$ and let $\cG
= I\cF_+I$. The operator $\cG$ is unitary and satisfies $\cG^2 = 1$. Let
$\Lambda>0$ and let $H_\Lambda = I(K_{1/\Lambda,1/\Lambda}) = L^2((0,\Lambda),
dt)\bigcap \cG(L^2((0,\Lambda), dt))$. The Mellin transforms of elements from
$H_\Lambda$ are entire functions with trivial zeros at $-2j$, $j\in\NN$. The
function $\wh{f}(s)$ satisfies the same functional equation as the Riemann zeta
function if and only if $\cG(f) = f$, and up to a sign if and only if $\cG(f) =
-f$. We prove:

\medskip

{\sc Theorem~2.\quad}{\it One has $\bigcap H_\Lambda = \{0\}$ and
$\overline{\bigcup H_\Lambda} = K$. For each $\Lambda>0$, each complex number
$w$, and each integer $k\in\NN$ there is a unique vector $Z^\Lambda_{w,k}$ in
$H_\Lambda$ such that the euclidean scalar product $(f,Z^\Lambda_{w,k}]$ equals
$(d^k/d^kw)\left(\pi^{-w/2}\Gamma({w/2})\wh{f}(w)\right)$ for each $f\in H_\Lambda$. The
vectors $Z^\Lambda_{w,k}$, $w\in\CC$, $k\in\NN$ are linearly independent, in
particular they are all non-vanishing.}

\medskip

Let $\Lambda>1$ and let $\cV_\Lambda$ be the vector space of smooth functions
with support in $[1/\Lambda, \Lambda]$. Let $D$ be the operator $ u (d^2/d^2u)
u$. Let, for $\phi\in\cV_\Lambda$, $E(\phi)(u)=\sum_{n\geq 1}\phi(nu) -
(\int_0^\infty\phi(t)dt)/u$.

\medskip {\sc Theorem~3.\quad}{\it Let $W_\Lambda$ be the closure in $K$ of the
functions $E(\phi)$ for $\phi\in D(\cV_\Lambda)$. One has $W_\Lambda\subset
H_\Lambda$. A vector $Z_{w,k}^\Lambda$ is perpendicular to $W_\Lambda$  if and
only if $w$ is a non-trivial zero $\rho$ of the zeta function of Riemann and
$0\leq k<m_\rho$ ($m_\rho$ = multiplicity of $\rho$.)}

\medskip

The outlines of the proofs will be found in the french part of this Note. One
thus defines $HP_\Lambda$ as the perpendicular complement  to $W_\Lambda$ in
$H_\Lambda$ so that $$K = \left(L^2((\Lambda,\infty),dt) + \cG
L^2((\Lambda,\infty),dt) \right) \perp W_\Lambda \perp HP_\Lambda$$  The
construction should be compared to the approach of Connes in [11], [12].

\medskip

\noindent{\bf Conclusion.\quad} In the continuation of [4], [5], we have shown
that the {\it scattering theory} for the Fourier transform (which belongs to
Analysis) has arithmetical aspects, through exhibiting some natural Hilbert
spaces containing vectors intrinsically associated to the non-trivial zeros of
the Riemann zeta function (and to their multiplicities; this is a continuation
to [10].) It would be interesting to better understand the operator framework
behind this construction ([6],[7],[8].) The study of $\Lambda \to 0$,
$\Lambda=1$, $\Lambda\to\infty$ is a problem of the {\it renormalization group}
([16]), which appears to be linked to the problem of the nature of the zeta
function ([3].)

\bigskip


\vskip 6pt
\centerline{\hbox to 2cm{\hrulefill}}
\vskip 9pt


Dans la section intitul\'ee ``Espaces de Sonine'' de son livre ([13]) 
De~Branges associe des espaces de Hilbert de fonctions enti\`eres \`a la
transformation de Hankel d'ordre $\nu$, $\nu>-1$. Lorsque $\nu=-1/2$ la
transformation de Hankel est essentiellement la transformation en cosinus
$\cF_+$ agissant sur $K = L^2(]0,\infty[,dt)$ selon (formellement)
$\cF_+(f)(t)=2\int_0^\infty \cos(2\pi tu)f(u)du$. Nous construisons dans la
pr\'esente Note des quotients de ces Espaces de Sonine (associ\'es \`a $\cF_+$)
qui contiennent des vecteurs intrins\`equement li\'es aux z\'eros non-triviaux
de la fonction dz\^eta de Riemann (et \`a leurs multiplicit\'es \'eventuelles.)
Une description plus satisfaisante de cette construction est de nature
op\'eratorielle, et nous reviendrons sur ce point dans un travail ult\'erieur.

\bigskip

\noindent{\bf 1. Les espaces $K_{a,b}$}

\medskip

Soit $K$ l'espace de Hilbert $L^2(]0,\infty[,dt)$ des fonctions \`a valeurs
complexes de carr\'es int\'egrables sur $]0,\infty[$. Nous noterons
$(\phi,\psi]$ pour le produit scalaire {\it euclidien} $\int \phi(u)\psi(u)du$. 
Soit $\cF_+$ la transformation en cosinus agissant formellement selon
$\cF_+(f)(t)=2\int_0^\infty \cos(2\pi tu)f(u)du$. Elle v\'erifie $\cF_+^2 = 1$,
$(\cF_+(\phi),\cF_+(\psi)] = (\phi,\psi]$. Soit $I$ agissant selon $I(f)(t) =
1/t\;f(1/t)$. Le compos\'e $\Gamma_+ = \cF_+\,I$ est invariant  et est donc
diagonalis\'e par l'isom\'etrie de Mellin-Plancherel $f\mapsto \wh{f}(s) = \int
f(t) t^{s-1}\,dt$ entre $K$ et l'espace des fonctions de carr\'es int\'egrables
sur $\Real(s) = 1/2$ pour la mesure $|ds|/2\pi$. Le {\sl multiplicateur
spectral} associ\'e est la fonction $\gamma_+(s) = 2(2\pi)^{-s}\cos(\pi
s/2)\Gamma(s)$. On a donc $\gamma_+(s)\gamma_+(1-s) = 1$ et, presque partout sur
la droite critique: $$\forall f\in K\qquad\wh{f}(s) =
\gamma_+(s)\;\wh{\cF_+(f)}(1-s)\eqno (1.1)$$

Soient $a$ et $b$ deux nombres r\'eels strictement positifs et consid\'erons le
sous-espace ferm\'e $K_{a,b}$ de $K$ des fonctions $f$ support\'ees par
$[a,\infty[$ et telles que $\cF_+(f)$ soit support\'ee par $[b,\infty[$. Il est
connu classiquement que $L^2(]0,a])+\cF_+(L^2(]0,b]))$ est un sous espace
ferm\'e propre de $K$ ({\it voir\/} [15], {\it sect. 2.9, p. 126-127\/}) et donc
que les espaces $K_{a,b}$ sont tous non-r\'eduits \`a $\{0\}$. Si $f$ est
presque partout nulle sur $]0,a]$ alors  $\int f(t) t^{s-1}\,dt$ est analytique
dans le demi-plan $\Real(s)<1/2$, et par (1.1) si $f\in K_{a,b}$ alors
$\wh{f}(s)$ existe en tant que fonction analytique dans $\Real(s)>1/2$. En fait:

\medskip

{\sc Th\'eor\`eme~1.1.\ {\rm (voir [13])}\quad}{\it Toute fonction d'un espace $K_{a,b}$ a une
transform\'ee de Mellin qui est une fonction enti\`ere de $s\in\CC$.}
\medskip

Il suffirait de  constater en effet que les fonctions
$\pi^{-(1-s)/2}\Gamma((1-s)/2)\wh{f}(s)$ forment (\`a un changement de variables
pr\`es) un des Espaces de Sonine de fonctions enti\`eres consid\'er\'es par
De~Branges dans [13]. Mais pour \'etablir en particulier le Th\'eor\`eme (1.5)
plus loin, nous ferons plut\^ot reposer la d\'emonstration sur l'identit\'e
\'el\'ementaire suivante, valable pour tout $f\in K_{a,b}$:
$$\int_0^\infty f(t)\,t^{s-1}\,dt = \int_b^\infty {(1-s)\int_a^\infty
\sin(2\pi\,ut)t^{s-2}\,dt - a^{s-1}{\sin(2\pi\,au)}\over \pi
u}\cF_+(f)(u)\,du\eqno(1.2)$$

Cette identit\'e est tout d'abord \'etablie pour
$\Real(s)<0$ et donne le prolongement analytique de $\wh{f}(s)$ au moins au
demi-plan $\Real(s)<1$. Puis (par exemple) (1.1) donne le prolongement
analytique au plan complexe.   L'\'equation (1.2) r\'esulte du calcul
(imm\'ediat) par int\'egration par parties:
$$\Real(w)<0\;\Rightarrow\;2\int_a^\infty \cos(2\pi\,ut)t^{w-1}\,dt =
{(1-w)\int_a^\infty \sin(2\pi\,ut)t^{w-2}\,dt - a^{w-1}{\sin(2\pi au)}\over \pi
u}$$

Avec $C_a(u,w) = 2\int_a^\infty \cos(2\pi\,ut)t^{w-1}\,dt$ et
$S_a(u,w) = 2\int_a^\infty \sin(2\pi\,ut)t^{w-1}\,dt$ on a la relation ci-dessus
entre $C_a(u,w)$ et $S_a(u,w-1)$ et aussi une relation entre $S_a(u,w)$ et
$C_a(u,w-1)$ qui \'etablissent que  ces fonctions sont des fonctions enti\`eres
de $w$ (pour chaque $u>0$ fix\'e) et aussi qu'elles sont $O(1/u)$ sur
$[1,\infty[$ uniform\'ement par rapport \`a $w\in\CC$ lorsque $|w|$ est born\'e.
L'usage que nous ferons de ce dernier point est dans la continuation de [10]
(motiv\'e par [2], {\it Lemme 6}.) Il va sans dire que ces calculs
simples (et ceux qui suivent) sont bien connus.

Soit $f(t)$ une fonction infiniment d\'erivable \`a support compact dans
$]0,\infty[$. Par Fubini on a $\int_a^\infty \cF_+(f)(t)\,t^{w-1}\,dt =
\int_0^\infty C_a(u,w)f(u)\,du$ pour $\Real(w)<0$ puis pour tout $w\in\CC$
puisque les deux membres sont des fonctions enti\`eres de $w$. Soit $D_a(u,w) =
2\int_0^a \cos(2\pi\,ut)t^{w-1}\,dt$  (pour $\Real(w)>0$.) Par Fubini on a
$\int_0^a \cF_+(f)(t)\,t^{w-1}\,dt = \int_0^\infty D_a(u,w)f(u)\,du$ et donc
$\wh{\cF_+(f)}(w) = \int_0^\infty (C_a(u,w)+D_a(u,w))f(u)\,du$ d'o\`u par (1.1)
pour $\Real(w)>0$:
$$C_a(u,w) = \gamma_+(w)u^{-w} - 2\int_0^a \cos(2\pi\,ut)t^{w-1}\,dt
= \gamma_+(w)u^{-w} -2\sum_{j=0}^\infty {(-1)^j\over (2j)!}(2\pi
u)^{2j}{a^{2j+w}\over 2j+w}\eqno(1.3)$$

Cette derni\`ere formule montre que $C_a(u,w)$ pour $w$ fix\'e est une fonction
analytique de $u$ dans $\CC \setminus ]-\infty,0]$,  au moins pour $w$ distinct
de $0$, $-2$, $-4$, \dots  En utilisant le fait que $\gamma_+(w)$ a un p\^ole
simple en $-2k$ ($k\in\NN$) avec un r\'esidu \'egal \`a
$2(2\pi)^{2k}(-1)^k(1/(2k)!)$, on voit de plus que $C_a(u,-2k)$ diff\`ere d'une
fonction enti\`ere de $u$ par un multiple non nul de $u^{2k}\log(u)$ et qu'elle
est donc \'egalement analytique en $u$ sur $\CC \setminus ]-\infty,0]$. On peut
affirmer, les seuls z\'eros de $\gamma_+(w)$ \'etant en $1$, $3$, $5$, \dots,
que les seules fonctions enti\`eres en $u$ parmi les  $C_a(u,w)$, $w\in\CC$,
sont obtenues pour $w\in\{\pm 1, \pm 3, \pm 5, \dots\}$ et que $C_a(u,w)$ est
dans $K$ si et seulement si soit $\Real(w)<{1/2}$ soit $w\in\{1, 3, 5,
\dots\}$. Nous avons ainsi \'etabli:

\medskip {\sc Lemme~1.2.\quad}{\it La fonction $C_a(u,w)$ pour $u>0$
fix\'e est une fonction enti\`ere de $w$.}

\medskip {\sc Lemme~1.3.\quad}{\it Pour $w\in\CC$ fix\'e $C_a(u,w)$ est
une fonction analytique de $u\in\CC \setminus ]-\infty,0]$. Elle est enti\`ere
si et seulement si $w\in\{\pm(1+2j), j\in \NN\}$. Elle est $O(1/u)$ sur
$[1,\infty[$, et cela uniform\'ement par rapport \`a $w$ lorsque $|w|$ est
born\'e. Pour $\Real(w)>0$, elle vaut $\gamma_+(w)\,u^{-w} - 2\int_0^a
\cos(2\pi\,ut)t^{w-1}\,dt$, et est donc $\gamma_+(w)\,u^{-w} + O(1)$ sur $]0,1]$
et cela uniform\'ement par rapport \`a $w$ pour $\Real(w)\geq\epsilon>0$. Elle
appartient \`a $K$ si et seulement si soit $\Real(w)<{1/2}$ soit $w\in\{1, 3, 5,
\dots\}$. On a $C_a(u,1+2j) = - \cF_+(\Un_{t\leq a}t^{2j})(u)$ pour $j\in\NN$.}

\medskip {\sc Th\'eor\`eme~1.4.\quad}{\it Soient $a>0$, $b>0$ fix\'es. Les
transform\'ees de Mellin $\wh{f}(s)$ des fonctions $f\in K_{a,b}$ poss\`edent en
$1$, $3$, $5$, \dots des {\it z\'eros triviaux} et ce sont leurs seuls z\'eros
communs.}

\medskip

{\sc D\'emonstration.\quad} Nous donnons ici une approche directe qui ne fait
pas appel aux r\'esultats g\'en\'eraux de [13]. Le fait que $\wh{f}(1+2j)=0$
($j\in\NN$) pour tout $f\in K_{a,b}$ se voit directement sur l'\'equation
(1.1.). Pour tout nombre complexe $w$ et toute fonction $f\in K_{a,b}$ on a
$\wh{f}(w) = (\cF_+(f), \phi_{a,b}^{w}] = (f,\cF_+(\phi_{a,b}^{w})]$ avec
$\phi_{a,b}^w(u) = \Un_{u\geq b}(u)\,C_a(u,w)$. Donc si $w$ est un z\'ero commun
\`a toutes les fonctions $\wh{f}$ alors  $\phi_{a,b}^w \in \cF_+(L^2(]0,a])) +
L^2(]0,b])$. Dans l'\'ecriture (unique) correspondante $\phi_{a,b}^w = f_{a,b} +
g_{a,b}$ la fonction $f_{a,b}(u)$ est une fonction enti\`ere de $u$  et elle est
\'egale pour $u>b$ \`a $C_a(u,w)$ qui est analytique sur $\CC \setminus
]-\infty,0]$. Elles sont donc \'egales et on en d\'eduit que $C_a(u,w)$
appartient \`a $\cF_+(L^2(]0,a]))$. En particulier elle est dans $K$ ce qui
implique par ce qui pr\'ec\`ede soit $w\in\{1+2j, j\in\NN\}$ soit
$\Real(w)<1/2$. Mais ce dernier cas est \`a exclure car alors $C_a(u,w)$
appartient \`a $\cF_+(L^2([a,\infty[))$ (et est non nulle).~$\bullet$

\medskip L'\'equation (1.1) montre que $1$, $3$, $5$, \dots sont simples en tant
que z\'eros communs aux fonctions $\wh{f}(s)$ pour $f\in K_{a,b}$. Il existe une
isom\'etrie canonique entre $K_{a, b}$ et $K_{\lambda,\lambda}$ pour $\lambda^2
= ab$ et nous supposerons dor\'enavant $a=b=\lambda>0$. Nous d\'efinissons des
fonctions $D_{w,k}(t)$ sur $]0,\infty[$ et des vecteurs $X_{w,k}^\lambda$ dans
$K$ pour $w\in\CC$ et $k\in\NN$ de la mani\`ere suivante:
$X_{w,k}^\lambda(t):=\Un_{t\geq\lambda}(t)D_{w,k}(t)$ avec, pour $\Real(w) >
1/2$, $D_{w,k}(t):=(\log(1/t))^{k}t^{-w}$ et, pour $\Real(w)\leq 1/2$,
$D_{w,k}(t):=(d^k/d^kw)C_\lambda(t,w)$. Pour tout $f\in K_{\lambda,\lambda}$ on
a pour $\Real(w)> 1/2$ $(f,X_{w,k}^\lambda] = (d^k/d^kw)\wh{f}(1-w)$ et pour
$\Real(w)\leq 1/2$ on a $(f,X_{w,k}^\lambda] =
(d^k/d^kw)(\gamma_+(w)\wh{f}(1-w))$. On notera en particulier la compensation
entre les p\^oles de $\gamma_+(w)$ et les z\'eros triviaux de $\wh{f}(1-w)$.

\medskip {\sc Th\'eor\`eme~1.5.\quad}{\it Soit $\lambda>0$. Les projections
orthogonales des vecteurs $X_{w,k}^\lambda$ sur $K_{\lambda,\lambda}$ sont
lin\'eairement ind\'ependantes.}

\medskip

{\sc D\'emonstration.\quad} Supposons qu'une combinaison lin\'eaire finie des
$X_{w,k}^\lambda$ soit dans $K_{\lambda,\lambda}^\perp = L^2(]0,\lambda]) +
\cF_+(L^2(]0,\lambda]))$. Comme dans la d\'emonstration pr\'ec\'edente, la
combinaison lin\'eaire correspondante des $D_{w,k}(t)$ (qui est analytique en
$t$ sur $]0,\infty[$) doit \^etre dans $\cF_+(L^2(]0,\lambda]))$. Elle ne peut
\^etre de carr\'e int\'egrable au voisinage de l'origine que si il n'y a aucune
contribution d'un couple $(k,w)$ avec $\Real(w) > 1/2$ car alors
$D_{w,k}(t)=(\log(1/t))^{k}t^{-w}$. On \'elimine pour la m\^eme raison les
couples $(w,k)$ v\'erifiant $\Real(w) = 1/2$ car alors  $D_{w,k}(t)$ a par le
Lemme (1.3) une singularit\'e dominante pour $t\to 0$ du type
$(\log(1/t))^{k}t^{-w}$. On remarque finalement pour $\Real(w) < 1/2$ par
l'\'equation (1.1) l'identit\'e $(\cF_+(f), X_{w,k}^\lambda] =
(d^k/d^kw)(\wh{f}(w))$ qui montre que $\cF_+(X_{w,k}^\lambda)$ et $(-1)^k
X_{1-w,k}$ ont la m\^eme projection sur $K_{\lambda,\lambda}$. On est donc
ramen\'e au cas pr\'ec\'edemment \'elimin\'e.~$\bullet$

\bigskip\eject

\noindent{\bf 2. Les espaces $H_\Lambda$ et les vecteurs $Z_{w,k}^\Lambda$}

\medskip Il est naturel pour les d\'eveloppements ult\'erieurs de travailler
plut\^ot avec les fonctions $\wh{f}(1-s)$, pour lesquelles les z\'eros triviaux
sont situ\'es en $0$, $-2$, $-4$, \dots. Soit $H_\Lambda$ (pour $\Lambda>0$)
l'espace des fonctions dans $K$ presque partout nulles pour $t>\Lambda$ et dont
les images sous $\cG := I\cF_+I$ sont presque partout nulles pour $t>\Lambda$.
On a $H_{\Lambda} = IK_{\lambda, \lambda}$ avec $\lambda=1/\Lambda$ et les
transform\'ees de Mellin sont les fonctions $\wh{f}(1-s)$, $f\in K_{\lambda,
\lambda}$. L'unitaire $\cG$ v\'erifie $\cG^2 = 1$ et laisse stable $H_\Lambda$
qui se d\'ecompose donc comme une somme perpendiculaire $H_\Lambda^+\perp
H_\Lambda^-$, les transform\'ees de Mellin des fonctions invariantes sous $\cG$
satisfaisant la m\^eme \'equation fonctionnelle que la fonction dz\^eta de
Riemann, tandis que celles dans $H_\Lambda^-$ la satisfont au signe pr\`es. Il
est clair que les espaces $H_\Lambda$ forment une cha\^{\i}ne croissante
d'espaces de Hilbert.

\medskip
{\sc Proposition~2.1.\quad}{\it On a $\bigcap H_\Lambda = \{0\}$ et
$\overline{\bigcup H_\Lambda} = K$.}

\medskip

La premi\`ere assertion est imm\'ediate, et pour la deuxi\`eme on peut par
exemple utiliser le fait que l'espace $\bigcup H_\Lambda$ est invariant sous les
translations multiplicatives $D_\theta: f(t)\mapsto
(1/\sqrt{\theta})f(t/\theta)$ ($0<\theta<\infty$), puisque $D_\theta(H_\Lambda)
\subset H_{\Lambda\theta}$ pour $\theta\geq1$ et $D_\theta(H_\Lambda) \subset
H_{\Lambda/\theta}$ pour $\theta\leq 1$. Or si $f(t)$ est une fonction
quelconque non-nulle de $H_1$ alors les z\'eros de $\wh{f}(s)$ sur la droite
critique forment un ensemble de mesure de Lebesgue nulle et par un th\'eor\`eme
de Wiener cela implique que les combinaisons lin\'eaires de ses translat\'ees
multiplicatives sont denses dans $K$.

Au vu des z\'eros triviaux en $0$, $-2$, \dots on associe \`a $f\in H_\Lambda$
non plus simplement $\wh{f}(s)$ mais la fonction enti\`ere
$M(f)(s):=\pi^{-s/2}\Gamma({s/2})\wh{f}(s)$. L'espace de Hilbert des fonctions
enti\`eres $M(f)(s)$ pour $f\in H_\Lambda$ est un Espace de Sonine au sens de
[13], [14] (ici, la droite critique joue le r\^ole de l'axe r\'eel dans [13],
[14]). On dispose des \'equations fonctionnelles $M(\cG(f))(s) = M(f)(1-s)$.

\medskip
{\sc Proposition~2.2.\quad}{\it Pour tout nombre complexe $w$ et pour tout
entier $k\geq0$ l'\'evaluation $f\mapsto M(f)^{(k)}(w)$ d\'efinit une forme
lin\'eaire continue sur $H_\Lambda$.}

\medskip
Cela d\'ecoule des calculs pr\'ec\'edant le Th\'eor\`eme (1.5) puisque ces
formes lin\'eaires s'expriment (selon un syst\`eme triangulaire) en fonction des
produits scalaires euclidiens avec les vecteurs $I(X_{w,k}^\lambda)$ ($\lambda =
1/\Lambda$). Il existe donc un unique vecteur $Z_{w,k}^\Lambda$ dans $H_\Lambda$
tel que $\forall f\in H_\Lambda\ M(f)^{(k)}(w) = (f,Z_{w,k}^\Lambda]$. Le {\it
noyau reproduisant} (analytique) $K^\Lambda(w,k;z,l) = (Z_{w,k}^\Lambda,
Z_{z,l}^\Lambda]$ est un objet important de l'Analyse pour lequel la th\'eorie
de De~Branges donne une formule exacte \`a partir d'une certaine fonction
$E_\Lambda(z)$ au sujet de laquelle nous souhaiterions pouvoir \^etre plus
explicite. Les vecteurs $Z_{w,k}^\Lambda$ sont reli\'es aux projections
orthogonales des vecteurs $I(X_{w,k}^\lambda)$ ($\lambda = 1/\Lambda$) sur
$H_\Lambda$ par un syst\`eme triangulaire inversible et donc par le Th\'eor\`eme
(1.5):

\medskip

{\sc Th\'eor\`eme~2.3.\quad}{\it Les vecteurs $Z_{w,k}^\Lambda$ pour $w\in\CC$,
$k\in\NN$ sont lin\'eairement ind\'ependants (en particulier ils sont tous non
nuls.)}

\bigskip

\noindent{\bf 3. Les espaces $W_\Lambda$ et $HP_\Lambda$ pour $\Lambda>1$}

\medskip

Soit $\Lambda>1$ et soit $\cV_\Lambda$ l'espace vectoriel des fonctions
infiniment diff\'erentiables et support\'ees par $[1/\Lambda, \Lambda]$. La
transform\'ee de Mellin $\wh{\phi}(s)$ d'une telle fonction $\phi$ est une
fonction enti\`ere. Le produit $\zeta(s)\wh{\phi}(s)$ est une fonction
m\'eromorphe, de carr\'e int\'egrable sur la droite critique et qui est la
transform\'ee de Mellin au sens $L^2$ de la fonction $E(\phi)(u)=\sum_{n\geq
1}\phi(nu) - (\int_0^\infty\phi(t)dt)/u$ ([5].) Soit $D = u (d^2/d^2u) u$
l'op\'erateur diff\'erentiel invariant de multiplicateur spectral $s(s-1)$. Les
fonctions de $\cV_\Lambda$ dans l'image de $D$ sont exactement celles
v\'erifiant $\wh{\phi}(0) = \wh{\phi}(1) = 0$. On d\'esignera par $W_\Lambda$
l'adh\'erence dans $K$ des fonctions $E(\phi)$ pour
$\phi\in D(\cV_\Lambda)$.

\medskip

{\sc Th\'eor\`eme~3.1.\quad}{\it Soit $\Lambda>1$. On a $W_\Lambda\subset
H_\Lambda$. Un vecteur $Z_{w,k}^\Lambda$ est perpendiculaire \`a $W_\Lambda$ si
et seulement si $w$ est un z\'ero non-trivial $\rho$ de la fonction dz\^eta de
Riemann et $0\leq k<m_\rho$ ($m_\rho$ = multiplicit\'e de $\rho$.)}

\medskip

{\sc D\'emonstration.\quad} On observe pour $\phi\in D(\cV_\Lambda)$ que
$E(\phi)$ a son support dans $]0,\Lambda]$. De plus $I(\phi)$ est aussi dans
$D(\cV_\Lambda)$ et donc $EI(\phi)$ a \'egalement son support dans
$]0,\Lambda]$. La transform\'ee de Mellin de $E(\phi)$ est
$\zeta(s)\wh{\phi}(s)$ et son image sous $\cF_+ I$ est
$\gamma_+(s)\zeta(s)\wh{\phi}(s)=\zeta(1-s)\wh{\phi}(s)$. Donc $\cG(E(\phi)) =
I\cF_+ I(E(\phi))$ a comme transform\'ee de Mellin $\zeta(s)\wh{\phi}(1-s)$
autrement dit $\cG(E(\phi))=E(I(\phi))$. Il existe des fonctions $\wh{\phi}(s)$
avec $\phi\in \cV_\Lambda$ prenant des valeurs et des d\'eriv\'ees quelconques
en un nombre fini de nombres complexes fix\'es \`a l'avance et avec cela le
Th\'eor\`eme est d\'emontr\'e.~$\bullet$

\medskip

{\sc D\'efinition~3.2.\quad}{\it Soit $\Lambda>1$. On d\'esigne par $HP_\Lambda$
le compl\'ement orthogonal de $W_\Lambda$ dans $H_\Lambda$.}

\medskip
On comparera \`a l'approche de Connes dans [11], [12]. On a donc
$$K = \left(L^2((\Lambda,\infty),dt) + \cG
L^2((\Lambda,\infty),dt) \right) \perp W_\Lambda \perp HP_\Lambda$$ 
On peut traiter de mani\`ere
analogue le cas d'une s\'erie
$L(s,\chi)$ de Dirichlet pour les caract\`eres primitifs et pairs. Il sera alors
inutile de se restreindre \`a $D(V_\Lambda)$ et il faudra tenir compte du
conducteur $q>1$. Pour une s\'erie associ\'ee \`a un caract\`ere impair la
construction utilisera non plus $\cF_+$ mais la transformation en sinus
$\cF_-$. 

\bigskip\goodbreak

\noindent{\bf 4. Conclusion} 

\medskip\hyphenation{scat-tering}

Dans la continuation de [4], [5], nous avons montr\'e que
la th\'eorie (appartenant \`a l'Analyse) de collision ({\it scattering\/}) pour
la transformation de Fourier a des aspects arithm\'etiques en associant de
mani\`ere intrins\`eque  un certain quotient $HP_\Lambda$ \`a la fonction
dz\^eta de Riemann (pour $\Lambda>1$) et en montrant qu'il est li\'e \`a ses
z\'eros non-triviaux (et \`a leurs multiplicit\'es \'eventuelles; ceci est dans
la continuation de [10].) Il serait int\'eressant de mieux comprendre ces
constructions en termes op\'eratoriels ([6],[7],[8].) L'\'etude de
$\Lambda\to0$, de $\Lambda\to\infty$, et de $\Lambda=1$  est un probl\`eme de
{\it groupe de renormalisation\/} ([16]) qui semble reli\'e \`a celui de la
nature de la fonction dz\^eta ([3].)

\vfill\eject


\noindent{}

\bigskip

\centerline{\bf R\'ef\'erences bibliographiques}
\medskip

{\eightpt\rm\baselineskip=14pt\parskip=0pt

\item{[1]}
        B\'aez-Duarte L., 
        {\it A class of invariant unitary operators\/},
         Adv. in Maths. {\bf 144} (1999), 1-12.

\item{[2]}
        B\'aez-Duarte L., Balazard M., Landreau B. et Saias E.,
        {\it Notes sur la fonction $ \zeta $ de Riemann~{\bf 3}\/},
         Adv. in Maths. {\bf 149} (2000), 130-144.

\item{[3]}
          Blake W.,
          {\it The Explicit Formula and a propagator\/},
          21 p., {\tt http://fr.arXiv.org/abs/math/9809119}

\newdimen\larg
\setbox0=\hbox{Burnol}
\larg=\wd0

\item{[4]}
         Burnol~J.-F.,
          {\it Scattering on the $p$-adic field and a trace formula\/},
          Int. Math. Res. Not. {\bf 2000:2} (2000), 57-70.

\item{[5]}
         \hbox to \larg{\leaders\hrule\hfill},
         {\it An adelic causality problem related to abelian L-functions\/},
        J. Number Theory {\bf 87} (2001), no. 2, 253-269.
        
\item{[6]}
         \hbox to \larg{\leaders\hrule\hfill},
         {\it The Explicit Formula and the conductor operator\/},
          28p., {\tt http://fr.arXiv.org/abs/math/9902080}

\item{[7]}
         \hbox to \larg{\leaders\hrule\hfill},
         {\it Sur les Formules Explicites I: analyse invariante\/},
         C. R. Acad. Sci. Paris {\bf 331} (2000), s\'erie I, 423-428.

\item{[8]}
          \hbox to \larg{\leaders\hrule\hfill},
         {\it Quaternionic gamma functions and their logarithmic derivatives as
spectral functions\/}, Math. Res. Lett. {\bf 8} (2001), no. 1-2, 209-224.
          
\item{[9]}
          \hbox to \larg{\leaders\hrule\hfill},
          {\it A note on Nyman's equivalent formulation for the Riemann
            Hypothesis}, 4p., Cont. Math. Ser.,  \`a para\^{\i}tre.

\item{[10]}
         \hbox to \larg{\leaders\hrule\hfill},
         {\it A lower bound in an approximation problem involving the zeros of
the Riemann zeta function\/}, 17p., mars 2001. Soumis \`a Adv. in Maths. {\tt
http://fr.arXiv.org/abs/math/0103058}

\item{[11]}
         Connes A., 
         {\it Formule de trace en g\'eom\'etrie non-commutative et hypoth\`ese
de Riemann}, C. R. Acad. Sci. Paris {\bf 323} (1996), s\'erie I, 1231-1236.

\setbox0=\hbox{Connes}
\larg=\wd0

\item{[12]}
         \hbox to \larg{\leaders\hrule\hfill}, 
         {\it Trace formula in non-commutative Geometry and
the zeros of the Riemann zeta function},
         Selecta Math. (N.S.) {\bf 5} (1999), no. 1, 29-106.

\item{[13]}
         De~Branges L., 
         {\it Espaces Hilbertiens de Fonctions Enti\`eres\/},
         Masson et Cie, Paris (1972)
         (trad. de {\it Hilbert spaces of entire functions}
         Prentice Hall Inc., Englewood Cliffs, 1968.)

\setbox0=\hbox{De~Branges}
\larg=\wd0

\item{[14]}
         \hbox to \larg{\leaders\hrule\hfill}, 
         {\it A conjecture which implies the Riemann hypothesis\/},
         J. Funct. Anal. {\bf 121} (1994), no. 1, 117-184.

\item{[15]}
         Dym H., McKean H.P.,
        {\it Fourier Series and Integrals\/},
         Academic Press, 1972.

\item{[16]}
         Wilson K.G.,
         {\it The renormalization group and critical phenomena},
         Rev. Mod. Phys. {\bf 55} (1983), 583-600.

}
\end